\documentclass[a4paper]{amsart}

\title[Strong generators of subregular W-algebra]{%
Strong generators of the subregular $\calW$-algebra $\calW^{K-N}(\frsl_N, f_{sub})$
and combinatorial description at critical level}
\author{Naoki Genra}
\author{Toshiro Kuwabara}

\usepackage{latexsym}
\usepackage{amsmath}
\usepackage{amsfonts}
\usepackage{amssymb}
\usepackage{amsxtra}
\usepackage{mathrsfs}
\usepackage [matrix,arrow]{xy}

\usepackage[colorlinks]{hyperref}

\newtheorem{definition}{Definition}[section]
\newtheorem{proposition}[definition]{Proposition}
\newtheorem{theorem}[definition]{Theorem}

\newtheorem{lemma}[definition]{Lemma}
\newtheorem{remark}[definition]{Remark}

\newcommand{\refdef}[1]{Definition~\ref{#1}}
\newcommand{\refprop}[1]{Proposition~\ref{#1}}
\newcommand{\refthm}[1]{Theorem~\ref{#1}}

\newcommand{\reflemma}[1]{Lemma~\ref{#1}}

\newcommand{\refeq}[1]{(\ref{#1})}

\newcommand{\refsec}[1]{Section~\ref{#1}}

\newcommand{\WalgCK}{\calW^{\K}_{\C[\K]}}
\newcommand{\WalgR}{\calW^{\K}_{\bR}}
\newcommand{\ValgCK}{\calV^{\K}_{\C[\K]}}
\newcommand{\ValgR}{\calV^{\K}_{\bR}}

\newcommand{\K}{K}
\newcommand{\Kv}{K_0}  
\newcommand{\CK}{\C[\K]}
\newcommand{\bR}{R}
\newcommand{\intQ}{\int\!Q}
\newcommand{\intA}{\int\!A}
\newcommand{\intY}{\int\!Y}
\newcommand{\intalpha}{\int\!\alpha}
\newcommand{\intbeta}{\int\!\beta}
\newcommand{\Wak}{\mu_{W}}
\newcommand{\FMS}{\mu_{\beta\gamma}}
\newcommand{\W}{W}
\newcommand{\gfW}{\widetilde{\W}^{(N)}}

\newcommand{\Ctwo}[1]{\overline{A}(#1)}
\newcommand{\Zhu}[1]{A(#1)}
\newcommand{\barH}{\overline{H}}
\newcommand{\barE}{\overline{E}}
\newcommand{\barF}{\overline{F}}
\newcommand{\barW}{\overline{\W}}
\newcommand{\bare}{\overline{e}}
\newcommand{\barmu}{\overline{\Mu}}

\newcommand{\R}[1]{\rho_{#1}} 
\newcommand{\F}[1]{\sigma_{#1}} 
\newcommand{\tilR}[1]{\sigma_{#1}(u)}
\newcommand{\Mu}{X} 
\newcommand{\OpX}{\xi}  
\newcommand{\EleX}{a} 
\newcommand{\XI}{\intY}

\newcommand{\C}{{\mathbb C}}

\newcommand{\Z}{{\mathbb Z}}

\newcommand{\bbS}{{\mathbb S}}

\newcommand{\bbV}{{\mathbb V}}

\newcommand{\calD}{\mathcal{D}}

\newcommand{\calH}{\mathcal{H}}

\newcommand{\calV}{\mathcal{V}}
\newcommand{\calW}{\mathcal{W}}

\newcommand{\frS}{\mathfrak{S}}

\newcommand{\frg}{\mathfrak{g}}
\newcommand{\frh}{\mathfrak{h}}

\newcommand{\frsl}{\mathfrak{sl}}

\DeclareMathOperator{\lEnd}{\operatorname{\mathscr{E}\kern-.1pc\mathit{nd}}}
\DeclareMathOperator{\gHom}{Hom}
\DeclareMathOperator{\lHom}{\operatorname{\mathscr{H}\kern-.1pc\mathit{om}}}

\DeclareMathOperator{\Dim}{dim}
\DeclareMathOperator{\Ker}{Ker}

\DeclareMathOperator{\Span}{Span}

\newcommand{\blkbar}{\raisebox{0.5ex}{\rule{2ex}{0.4pt}}}

\newcommand{\One}{\mathbf{1}}
\newcommand{\NO}{{\genfrac{}{}{0pt}{1}{\circ}{\circ}}} 

\newcommand{\isoto}[1][]{\mathop{\xrightarrow[#1]%
{\rule{0pt}{.9ex}%
{\raisebox{-.4ex}[0ex][-.6ex]{$\mspace{3mu}\sim\mspace{3mu}$}}}}}
\newcommand{\scbul}{{\,\raise1pt\hbox{$\scriptscriptstyle\bullet$}\,}}

\begin{document}
\maketitle

\markboth{N.Genra and T.Kuwabara}{Strong generators of the subregular $\calW$-algebra}

\begin{abstract}
 We construct explicitly strong generators of the affine $\calW$-algebra 
 $\calW^{\Kv-N}(\frsl_N, f_{sub})$ of subregular type $A$. 
 Moreover, we are able to describe the OPEs between them at critical level.
 We also give a description the affine $\calW$-algebra $\calW^{-N}(\frsl_N, f_{sub})$
 in terms of certain fermionic fields, which was conjectured by Adamovi\'{c}. 
\end{abstract}

\section{Introduction}
\label{sec:intro}

For a reductive Lie algebra $\frg$, a nilpotent element $f \in \frg$ and
$k \in \C$, the affine $\calW$-algebra $\calW^{k}(\frg, f)$ is defined as
a vertex algebra constructed by the generalized quantum Drinfeld-Sokolov
reduction; see \cite{FF}, \cite{KRW}, \cite{KW1}. In this paper, we discuss
the affine $\calW$-algebra $\calW^{\Kv-N}(\frsl_N, f_{sub})$ associated with
$\frsl_N$ and a subregular nilpotent element $f_{sub} \in \frsl_N$ with
level $\Kv-N$, which we call the subregular $\calW$-algebra.
Recently, in \cite[Section 6]{Genra1}, the first author described the subregular 
$\calW$-algebra by using certain screening operators, and showed that the 
subregular $\calW$-algebra is isomorphic to a vertex algebra $W^{(2)}_N$ 
introduced by Feigin and Semikhatov in \cite{FS}.

For a principal nilpotent element $f_{pr} \in \frsl_N$, the corresponding affine 
$\calW$-algebra
$\calW^{\Kv-N}(\frsl_N, f_{pr})$ is a vertex algebra such that, at critical level $\Kv=0$,
$\calW^{-N}(\frsl_N, f_{pr})$ coincides with the center of affine vertex algebra $V^{-N}(\frsl_N)$, called the
Feigin-Frenkel center.
In \cite[Section 2]{AM}, Arakawa and Molev explicitly constructed strong 
generators of the vertex algebra $\calW^{\Kv-N}(\frsl_N, f_{pr})$. Their images 
through the Miura map are described by a certain noncommutative analog of 
the elementary symmetric polynomials, which recovers a result of Fateev and Lukyanov in \cite{FL}.

In \refsec{sec:strong-gen} of this paper, we discuss construction of certain strong generators for 
the subregular $\calW$-algebra $\calW^{\Kv-N}(\frsl_N, f_{sub})$. Our construction is based on 
the Feigin-Semikhatov description of $\calW^{\Kv-N}(\frsl_N, f_{sub})$, which describes
$\calW^{\Kv-N}(\frsl_N, f_{sub})$ as intersection of the kernels of screening
operators on a certain lattice vertex algebra.
We construct elements $\W_m$ ($m=2$, $\dots$, $N$) of the lattice vertex algebra by using 
the noncommutative elementary symmetric polynomials of the elements of the Heisenberg part 
of the lattice vertex algebra,
and show that they lie in the intersection of the kernels of the screening operators (\refdef{def:generators}
and \refprop{prop:W-belong}). 
We also show that these elements are algebraically independent in the Zhu's $C_2$ Poisson
algebra of the vertex algebra $\calW^{\Kv-N}(\frsl_N, f_{sub})$, and it implies that the
elements $\W_2$, $\dots$, $\W_{N-1}$, together with the generators $E$, $H$, $F$ of Feigin-Semikhatov's, 
strongly generate $\calW^{\Kv-N}(\frsl_N, f_{sub})$ (\refthm{thm:strong-gen}).

In \refsec{sec:critical-level}, we discuss the subregular $\calW$-algebra $\calW^{-N}(\frsl_N, f_{sub})$
at critical level, $\Kv = 0$. At critical level, the vertex algebra $\calW^{-N}(\frsl_N, f_{sub})$ has a nontrivial
center as a vertex algebra, and the center is naturally isomorphic to the Feigin-Frenkel center of
the affine vertex algebra $V^{-N}(\frsl_N)$.
We show that our elements $\W_2$, $\dots$, $\W_N$ strongly generate the center (\refprop{prop:gen-center}).
Moreover, we give an explicit form of OPEs between the strong generators (\refthm{thm:OPE}).

In \cite{Adamovic}, Adamovi\'{c} conjectured that the subregular $\calW$-algebra 
$\calW^{-N}(\frsl_N, f_{sub})$ is isomorphic to a vertex algebra generated by certain
fields, consisting of certain fermionic fields and the generators of the Feigin-Frenkel center.
We prove his conjecture by using our strong generators (\refthm{thm:Adamovic}).

\subsection*{Acknowledgments}

One of the main results of this paper, \refthm{thm:Adamovic}, was first conjectured by Dra\v{z}en Adamovi\'{c}.
The authors would like to express their gratitude to him for showing us his private notes \cite{Adamovic}.
The authors are deeply grateful to Tomoyuki Arakawa for valuable comments.
The authors also thank Boris Feigin and Alexei Semikhatov for fruitful discussion on their construction
of the vertex algebra $W^{(2)}_N$. 
The second author thanks Yoshihiro Takeyama for discussion on the proof of \reflemma{lemma:2}.

The first author was supported by Grant-in-Aid for JSPS Fellows (No.17J07495). 
The second author was supported by JSPS KAKENHI Grant Number JP17K14151.

\section{Subregular $\calW$-algebras and Feigin-Semikhatov screenings}
\label{sec:def-subreg-W}
The subregular $\calW$-algebra $\calW^{\Kv}=\calW^{\Kv-N}(\frsl_N, f_{sub})$ at level $\Kv-N$
is a vertex algebra defined by 
the generalized quantum Drinfeld-Sokolov reduction associated with $\frsl_N, f_{sub}$ and $\Kv\in\C$ \cite{KRW}, 
where $f_{sub} = e_{- \alpha_2} + \dots + e_{- \alpha_{N-1}} \in \frsl_N$ is a subregular nilpotent element in 
$\frsl_N$. We introduce a free field realization of $\calW^{\Kv}$ following \cite{FS} and \cite{Genra1}.

We follow \cite{FB, Kac} for definitions of vertex algebras, and denote by $A(z)=Y(A,z)=\sum_{n\in\Z} A_{(n)}$
a field on $V$ for a element $A$ in a vertex algebra $V$. Let $\K$ be an indeterminate and $\WalgCK$
the subregular $\calW$-algebra associated with $\frsl_N, f_{sub}, \K$ over $\CK$. By definition, we have 
$\WalgCK \otimes \C_{\Kv}=\calW^{\Kv}$, where $\C_{\Kv} = \C[\K]/(\K-\Kv) \simeq \C$ is a one-dimensional 
$\C[\K]$-module on which $\K$ acts by $\Kv$. See e.g. \cite{ACL2}. Let 
$\bbV = \CK A_{N-1} \oplus \dots \oplus \CK A_{1} \oplus \CK Q \oplus \CK Y$ be
a free $\CK$-module of rank $N+1$.
We define a symmetric bilinear form on $\bbV$ given by the Gram matrix
\[
 \begin{pmatrix}
  2\K & -\K & 0 & 0 & \cdots & \cdots & \cdots & \cdots & 0 \\
  -\K & 2\K & -\K & 0 & \cdots & \cdots & \cdots & \cdots & 0 \\
  0 & -\K & 2\K & -\K & 0 & \cdots & \cdots & \cdots &  0 \\
  \cdots & \cdots & \cdots & \cdots & \cdots & \cdots & \cdots & \cdots & \cdots \\
  \cdots & \cdots & \cdots & \cdots & \cdots & \cdots & \cdots & \cdots & \cdots \\
  \cdots & \cdots & \cdots & \cdots & \cdots & \cdots & \cdots & \cdots & \cdots \\
  0 & \cdots & \cdots & \cdots & 0 & -\K & 2\K & -\K & 0\\
  0 & \cdots & \cdots & \cdots & \cdots & 0 & -\K & 1 & 1 \\
  0 & \cdots & \cdots & \cdots & \cdots & 0 & 0 & 1 & 0 \\
 \end{pmatrix}.
\]
Consider a Heisenberg vertex algebra $\calH^{\K}$ over $\CK$ associated with the bilinear
form on $\bbV$. It is a vertex algebra
generated by the elements $A_i$, $Q$ and $Y$ ($i=1$, $\dots$, $N-1$) subject to the OPE
$a(z) b(w) \sim (a, b) / (z-w)^2$ where $a$, $b = A_i$, $Q$ or $Y$, and $(\;\,, \;)$
is the bilinear form on $\bbV$. For a vector $w \in \bbV$,
let $\calH^{\K}_{w}$ be the $\calH^{\K}$-module of highest weight
$w \in \bbV$. We denote the anti-derivation of $w$ by $\int\!w$
and the highest weight vector of $\calH^{\K}_{w}$ by $e^{\int\!w}$. The direct sum
$\ValgCK = \bigoplus_{m \in \Z} \calH^{\K}_{m Y}$ is equipped with a vertex algebra
structure. Indeed, the vertex operator $e^{m \XI}(z) = \sum_{n} e^{m \XI}_{(n)} z^{-n-1}$ corresponding to 
the highest weight vector $e^{m \XI}$ is a field with the following OPEs
\[
 e^{m \XI}(z) e^{n \XI}(w) \sim 0, \quad a(z) e^{m \XI}(w) \sim \frac{(a, mY)}{z-w}\, e^{m \XI}
\]
for $m$, $n \in \Z$, and $a \in \bbV$ and the derivative of $e^{m \XI}(z)$ is given by
$\partial e^{m \XI}(z) = \NO m Y(z) e^{m \XI}(z) \NO$.

Similarly to $e^{\XI}(z)$, we also have the vertex operator $e^{\intQ}(z)$ (resp. $e^{\intA_i}(z)$)
associated with $Q \in \bbV$ (resp. $A_i \in \bbV$ for $i=1$, $\dots$, $N-1$). The OPEs
between these vertex operators and fields of $\ValgCK$ are given as follows:
\begin{align*}
&a(z) e^{\int b}(w) \sim \frac{(a, b)}{z-w} e^{\int b}(w), \quad
e^{\intA_i}(z) e^{m \XI}(w) \sim 0, \\
&e^{\intQ}(z) e^{m \XI}(w) \sim (z-w)^{m} e^{\pm m \XI + Q}(w)
\end{align*}
where $a$, $b \in \bbV$, $i=1$, $\dots$, $N-1$ and $m \in \Z$.
The residue
of $e^{\intQ}(z)$ (resp. $e^{\intA_i}$(z)) gives an operator on $\ValgCK$ such that
$e^{\intQ}_{(0)} : \calH^{\K}_{m Y} \longrightarrow \calH^{\K}_{m Y + Q}$
(resp. $e^{\intA_i}_{(0)}: \calH^{\K}_{m Y} \longrightarrow \calH^{\K}_{m Y + A_i}$) for
$m \in \Z$ and $i=1$, $\dots$, $N-1$. The operators $e^{\intQ}_{(0)}$, $e^{\intA_i}_{(0)}$
are called screening operators. 
The vertex algebra given as intersection of the kernels of these screening operators were introduced by
Feigin and Semikhatov in \cite{FS}. Recently, the first author showed that their vertex algebra
is isomorphic to the subregular $\calW$-algebra.

\begin{proposition}[\cite{Genra1}, Theorem 6.9]
\label{prop:G-isom}
 As a vertex algebra over the ring $\CK$, we have an isomorphism
\[
 \mu^{\K} \colon \WalgCK \isoto \Ker e^{\intQ}_{(0)} \cap \bigcap_{i=1}^{N-1} \Ker e^{\intA_i}_{(0)}.
\]
\end{proposition}

Since $\calW^{\K}_{\CK} \otimes \C_{\Kv}=\calW^{\Kv}$, we have an embedding 
$\calW^{\Kv} \hookrightarrow \calV^{\Kv}$ for all $\Kv\in\C$. We remark that the embedding is obtained as 
the composition of three maps $\mu, \Wak, \FMS$ defined as follows. Applying the specialization functor 
$?\otimes \C_{\Kv}$ to embeddings
\begin{align*}
\Ker e^{\intQ}_{(0)} \cap \bigcap_{i=1}^{N-1} \Ker e^{\intA_i}_{(0)} \rightarrow \Ker e^{\intQ}_{(0)} \cap 
\Ker e^{\intA_1}_{(0)} \rightarrow \Ker e^{\intQ}_{(0)} \rightarrow \ValgCK,
\end{align*}
we have vertex algebra homomorphisms
\begin{multline*}
 \calW^{\Kv} \xrightarrow{\,\,\mu\,\,} V^{\tau_{\Kv-N}}(\frg_0) \simeq V^{\Kv-2}(\frsl_2) \otimes V^{\Kv}(\C^{N-2}) \\
 \xrightarrow{\,\Wak\,} \calD^{ch}(\C^1) \otimes V^{\Kv}(\C^{N-1})
 \xrightarrow{\,\FMS\,} \calV^{\Kv},
\end{multline*}
where $V^{\tau_{\Kv-N}}(\frg_0)$ is the affine vertex algebra associated with
the Lie subalgebra $\frg_0$ and the bilinear form $\tau_{\Kv-N}$ on $\frg_0$ defined
in \cite[(2.2)]{Genra1}, and $\calD^{ch}(\C^1)$ is the vertex algebra of $\beta\gamma$-system of rank one.
It then follows that $\mu, \Wak, \FMS$ are injective maps, called the Miura map for $\calW^{\Kv}$
\cite{KW1, Genra1}, Wakimoto realization for $V^{\Kv-2}(\frsl_2)$ \cite{Wakimoto, Frenkel} and Friedan-Martinec-Shenker
bosonization \cite{FMS} respectively.

\section{Strong generators of the subregular $\calW$-algebra}
\label{sec:strong-gen}

In this section, we explicitly construct elements $\W_1$, $\dots$, $\W_N$ of the subregular
$\calW$-algebra, and show that $E$, $F$, $H$, $\W_2$, $\dots$, $\W_{N-1}$ strongly generate
the vertex algebra $\WalgCK \otimes \C_{\Kv}$ for any $\Kv \in \C \setminus \{1\}$.
In the rest of this paper, we extend the vertex algebra $\WalgCK$ (resp. $\ValgCK$) by the ring 
$\bR := \C[\K, (\K-1)^{-1}]$, and write $\WalgR = \WalgCK \otimes_{\C[\K]} R$
(resp. $\ValgR = \ValgCK \otimes_{\C[\K]} R$) for short.

Set $\ell_N(\K) = \K (N-1) / N - 1$. Define three elements in $\calV^{\K}$
\begin{gather*}
 H = \ell_N(\K) Y + Q + \frac{N-1}{N} A_1 + \dots + \frac{2}{N} A_{N-2} + \frac{1}{N} A_{N-1}, \\
 E = e^{\XI}, \qquad F = - \R{N} \dots \R{2} \R{1} e^{- \XI}
\end{gather*}
where $\R{i} = (\K-1)(\partial + Y_{(-1)}) + Q_{(-1)} + \sum_{j=1}^{i-1} A_{j (-1)}$ for $i=1$, $\dots$, $N$.
Note that $\R{1} e^{- \XI} = Q_{(-1)} e^{- \XI}$ because $\partial e^{- \XI} = - Y_{(-1)} e^{- \XI}$. 
In \cite{FS}, they showed that these three elements generate the vertex algebra $\WalgR$. 
The first goal of the present paper is to construct a set of strong generators of $\WalgR$,
including these three elements $H$, $E$ and $F$.

Define $N$ elements in the Heisenberg part $\calH^{\K}$ of the vertex algebra $\ValgR$
\[
 \Mu_i = - \frac{\K}{N} Y - \sum_{j=1}^{i-1} \frac{j}{N} A_j + \sum_{j=i}^{N-1} \frac{N-j}{N} A_j \in \calH^{\K}
\]
for $i=1$, $\dots$, $N$. Then, we have $\R{i} = (\K-1) \partial + H_{(-1)} - \Mu_{i (-1)}$ for $i=1$, $\dots$, $N$.
Also, we define 
\[
 \Mu_0 = - \frac{\K}{N} Y + Q + \frac{N-1}{N} A_1 + \dots + \frac{1}{N} A_{N-1} \in \calH^{\K} \subset \ValgR
\]
and $\R{0} = (\K-1) \partial + H_{(-1)} - \Mu_{0 (-1)} = (\K-1) (\partial + Y_{(-1)})$.

Recall the definition of the noncommutative elementary symmetric polynomials (cf. \cite[(12.48)]{Molev}).
Let $\OpX_1$, $\dots$, $\OpX_N$ be mutually noncommutative $N$ operators on a certain vector space. 
Define the $m$-th noncommutative elementary symmetric polynomial in $\OpX_1$, $\dots$, $\OpX_N$,
\[
 e_m(\OpX_1, \dots, \OpX_n) = \sum_{i_1 > i_2 > \dots > i_m} \OpX_{i_1} \OpX_{i_2} \dots \OpX_{i_m}.
\]
Note that we arrange the operators reverse-lexicographically. 

Let $\F{i} = (\K-1)\partial - \Mu_{i (-1)}$ be an operator on $\ValgR$ for $i=0$, $1$, $\dots$, $N$.
For $m=1$, $\dots$, $N$, we define an element of $\ValgR$ 
\begin{equation}
\label{eq:def-Wm1}
 \W'_m = \sum_{k=0}^{N} (-1)^{k} \prod_{j=1}^{k} \frac{j(\K-1) + 1}{j(\K-1)}
 \hspace*{-3mm} \sum_{1 \le i_1 < \dots < i_{N-k} \le N} \hspace*{-8mm}
 e_m(\overbrace{\F{0}, \dots, \F{0}}^{\text{$k$-times}}, \F{i_1}, \dots, 
 \F{i_{N-k}}). 
\end{equation}

We also introduce the generating function of these elements $\W'_m$ ($m=1$, $\dots$, $N$).
Let $u$ be an indeterminate which commutes with all other elements. Note that, for operators 
$\OpX_1$, $\dots$, $\OpX_N$, we have
\[
 (u + \OpX_N) \dots (u + \OpX_1) 
 = \sum_{m=0}^{N} e_m(\OpX_1, \dots, \OpX_N) \, u^{N-m}
\]
by the definition of $e_m$.
Setting $\tilR{i} = u + \F{i}$, we have
\begin{multline}
 \gfW(u) := \sum_{m=0}^{N} \W'_m u^{N-m} \\
 = \sum_{k=0}^{N} (-1)^{k} \prod_{j=1}^{k} \frac{j(\K-1) + 1}{j(\K-1)} \hspace*{-3mm}
 \sum_{N \ge i_1 > \dots > i_{N-k} \ge 1} \hspace*{-8mm} \tilR{i_1} \dots \tilR{i_{N-k}} \tilR{0}^{k} \One
\end{multline}
where $\W'_0$ is equal to $\One$ up to multiplication by a certain constant.

We will show that the elements $\W'_m$ ($m=1$, $\dots$, $N$) belong to the subregular $\calW$-algebra
$\WalgR$ by calculating the action of the screening operators $e^{\intQ}_{(0)}$, 
$e^{\intA_i}_{(0)}$ ($i=1$, $\dots$, $N-1$). 

\begin{lemma}
 \label{lemma:comm-rel}
 (1) For $i=1$, $\dots$, $N-1$, $j=0$, $\dots$, $N$ and $m \in \Z$, we have
 \[
  [\,e^{\intA_i}_{(m)}, \tilR{j}\,] = 
 \begin{cases}
 (m (\K-1) + \K) e^{\intA_i}_{(m-1)} & (j=i) \\
 (m (\K-1) - \K) e^{\intA_i}_{(m-1)} & (j=i+1) \\
 m (\K-1) e^{\intA_i}_{(m-1)} & (j \ne i, i+1)
 \end{cases}
 \]
 in $\bigoplus_{m \in \Z} \gHom_{\C[\K]}(\calH^{\K}_{m Y}, \calH^{\K}_{m Y + A_i})[u]$.

 (2) For $j=0$, $\dots$, $N$ and $m \in \Z$, we have
 \[
  [\,e^{\intQ}_{(m)}, \tilR{j}\,] = 
 \begin{cases}
 (m+1) (\K-1) e^{\intQ}_{(m-1)} & (j=0) \\
 (m (\K-1) - \K) e^{\intQ}_{(m-1)} & (j=1) \\
 m (\K-1) e^{\intQ}_{(m-1)} & (j = 2, \dots, N)
 \end{cases}
 \]
  in $\bigoplus_{m \in \Z} \gHom_{\C[\K]}(\calH^{\K}_{m Y}, \calH^{\K}_{m Y + Q})[u]$.
\begin{proof}
 Note that $(A_i, \Mu_i) = \K$, $(A_i, \Mu_{i+1}) = - \K$, $(A_i, \Mu_j) = 0$ for $j \ne i, i+1$,
 and $(Q, \Mu_0) = - \K + 1$, $(Q, \Mu_1) = - \K$, $(Q, \Mu_j) = 0$ for $j = 2$, $\dots$, $N$.
 Then, both (1) and (2) can be checked by direct computation.
\end{proof}
\end{lemma}

\begin{proposition}
 \label{prop:1}
 For $i=1$, $\dots$, $N-1$, we have $e^{\intA_i}_{(0)} \gfW(u) = 0$.
\begin{proof}
 First, note that we have 
\[
 e^{\intA_i}_{(0)} \tilR{N} \dots \tilR{i+2} \tilR{i+1} \dots \One = \tilR{N} \dots \tilR{i+2} e^{\intA_i}_{(0)} \tilR{i+1} \dots \One
\]
 since the screening operator $e^{\intA_i}_{(0)}$ commutes with $\tilR{j}$ for 
 $j \ne i$, $i+1$ by \reflemma{lemma:comm-rel}.

 In $\gfW(u)$, there exists three kinds of terms; i) $\cdots \tilR{i+1} \tilR{i} \cdots \One$, terms with both
 factors $\tilR{i+1}$ and $\tilR{i}$, ii) $\cdots \tilR{i+1} \cdots \One$ or $\cdots \tilR{i} \cdots \One$,
 terms with either $\tilR{i+1}$ or $\tilR{i}$, iii) $\cdots (\tilR{i+1})^{\wedge} (\tilR{i})^{\wedge} \cdots \One$,
 terms without $\tilR{i+1}$ nor $\tilR{i}$. We consider the action of $e^{\intA_i}_{(0)}$ in these
 three cases individually.

 i) By \reflemma{lemma:comm-rel} (1), we have
 \begin{multline*}
  e^{\intA_i}_{(0)} \cdots \tilR{i+1} \tilR{i} \cdots \One \\
  = \cdots [\, e^{\intA_i}_{(0)}, \tilR{i+1} \,] \tilR{i} \cdots \One +
  \cdots \tilR{i+1} [\, e^{\intA_i}_{(0)}, \tilR{i} \,] \cdots \One \\
  \shoveleft{= \K \cdots (\tilR{i+1} - \tilR{i}) e^{\intA_i}_{(-1)} \cdots \One} \\
  - \K \cdots ((-1)(\K-1) + \K) e^{\intA_i}_{(-2)} \cdots \One \\
  = \K \cdots \bigl\{ A_{i (-1)} e^{\intA_i}_{(-1)} - e^{\intA_i}_{(-2)} \bigr\} \cdots \One = 0.
 \end{multline*}
 Here we used $[e^{\intA_i}_{(m)}, \Mu_{j (n)}] = 0$ for $j \ne i$, $i+1$ with $m$, $n \in \Z$,
 and $e^{\intA_i}_{(-2)} = \bigl(A_{i (-1)} e^{\intA_i} \bigr)_{(-1)}$.

 ii) Note that the generating function $\gfW(u)$ is symmetric by permutation between factors $\tilR{1}$,
 $\dots$, $\tilR{N}$. This implies that, for a term $(\cdots)_1 \tilR{i+1} (\cdots)_2 \One$ of type ii),
 we also have a term $(\cdots)_1 \tilR{i} (\cdots)_2 \One$ of exactly the same form except for replacing
 $\tilR{i+1}$ by $\tilR{i}$. Here the factors different from $\tilR{i+1}$ and $\tilR{i}$ are denoted by
 $(\cdots)_j$ ($j=1$, $2$). Then, by \reflemma{lemma:comm-rel} (1), we have
\begin{multline*}
  e^{\intA_i}_{(0)} \bigl\{(\cdots)_1 \tilR{i+1} (\cdots)_2 \One + (\cdots)_1 \tilR{i} (\cdots)_2 \One \bigr\} \\
 = (\cdots)_1 \bigl(- \K e^{\intA_i}_{(-1)} \bigr) (\cdots)_2 \One
 + (\cdots)_1 \bigl(+ \K e^{\intA_i}_{(-1)} \bigr) (\cdots)_2 \One = 0.
\end{multline*}
 
 iii) A term without $\tilR{i+1}$ and $\tilR{i}$ trivially vanishes by the action of the screening operator
 $e^{\intA_i}_{(0)}$ by \reflemma{lemma:comm-rel} (1).

 As a consequence, we have $e^{\intA_i}_{(0)} \gfW(u) = 0$.
\end{proof}
\end{proposition}

The action of another screening operator $e^{\intQ}_{(0)}$ is more complicated. Previous to 
the calculation of $e^{\intQ}_{(0)} \gfW(u)$  we prepare the following lemma.

\begin{lemma}
 \label{lemma:1}
 For $m \ge 0$, we have
 \[
  e^{\intQ}_{(0)} \Bigl( \tilR{1} \tilR{0}^m \One - \frac{(m+1)(\K-1) + 1}{(m+1)(\K-1)} \tilR{0}^{m+1} \One \Bigr) = 0
 \]
 \begin{proof}
  First, note that 
  \begin{align*}
   [\, e^{\intQ}_{(0)}, \tilR{0}^l \,] &= \sum_{i=1}^{l} (-1)^i (\K-1)^{i} \frac{l!}{(l-i)!} 
   \tilR{0}^{l-i} e^{\intQ}_{(-i)}, \\
   [\, e^{\intQ}_{(-1)}, \tilR{0}^l \,] &= \sum_{i=2}^{l+1} (-1)^{i-1} (\K-1)^{i-1} \frac{l!}{(l-i+1)!} 
   i \tilR{0}^{l-i+1} e^{\intQ}_{(-i)},
  \end{align*}
  by \reflemma{lemma:comm-rel} (2).

  First we deal with the first term of the equality of the lemma.
  Using the fact $\tilR{1} = \tilR{0} + Q_{(-1)}$ and the identity
  \[
  Q_{(-1)} \tilR{0}^{m-i} = \sum_{j=0}^{m-i} (-1)^j \frac{(m-i)!}{(m-i-j)!} \tilR{0}^{m-i-j} Q_{(-j-1)},
  \]
  we obtain 
  \begin{multline}
   \label{eq:1}
   e^{\intQ}_{(0)} \tilR{1} \tilR{0}^m \One = 
   \tilR{1} e^{\intQ}_{(0)} \tilR{0}^m \One - \K e^{\intQ}_{(-1)} \tilR{0}^m \One \\
   \shoveleft{= - \K \tilR{0}^m e^{\intQ}_{(-1)} \One} \\
   + \K \sum_{i=2}^{m+1} (-1)^i (\K-1)^{i-1} \frac{m!}{(m-i+1)!} i \tilR{0}^{m-i+1} e^{\intQ}_{(-i)} \One \\
   + \sum_{i=1}^{m} (-1)^i (\K-1)^{i} \frac{m!}{(m-i)!} \tilR{1} \tilR{0}^{m-i} e^{\intQ}_{(-i)} \One \\
   \shoveleft{= \sum_{i=1}^{m+1} (-1)^i (\K-1)^{i-1} \frac{m!}{(m-i+1)!} i \tilR{0}^{m-i+1} e^{\intQ}_{(-i)} \One} \\
   + \sum_{i=1}^{m} (-1)^i (\K-1)^{i} \frac{(m+1)!}{(m-i+1)!} \tilR{0}^{m-i+1} e^{\intQ}_{(-i)} \One \\
   + \sum_{i=1}^{m} (-1)^i (\K-1)^{i} \frac{m!}{(m-i)!} Q_{(-1)} \tilR{0}^{m-i} e^{\intQ}_{(-i)} \One \\
   \shoveleft{= \sum_{i=1}^{m+1} (-1)^i (\K-1)^{i-1} \frac{m!}{(m-i+1)!} i \tilR{0}^{m-i+1} e^{\intQ}_{(-i)} \One} \\
   + \sum_{i=1}^{m} (-1)^i (\K-1)^{i} \frac{(m+1)!}{(m-i+1)!} \tilR{0}^{m-i+1} e^{\intQ}_{(-i)} \One \\
   + \sum_{i=1}^{m} \sum_{j=0}^{m-i} 
   (-1)^{i+j} (\K-1)^{i+j} \frac{m!}{(m-i-j)!} \tilR{0}^{m-i-j} Q_{(-j-1)} e^{\intQ}_{(-i)} \One.
  \end{multline}
  
  On the other hand, for the second term of the equality of the lemma, we have
  \begin{multline}
   \label{eq:2}
   e^{\intQ}_{(0)} \frac{(m+1)(\K-1)+1}{(m+1)(\K-1)} \tilR{0}^{m+1} \One \\
   = \sum_{i=1}^{m+1} (-1)^i (\K-1)^i \frac{(m+1)!}{(m-i+1)!} \tilR{0}^{m+1-i} e^{\intQ}_{(-i)} \One \\
   + \sum_{i=1}^{m+1} (-1)^i (\K-1)^{i-1} \frac{m!}{(m-i+1)!} \tilR{0}^{m+1-i} e^{\intQ}_{(-i)} \One \\
  \end{multline}
  
  The second term in \refeq{eq:1} and the first term in \refeq{eq:2} are canceled out. 
  Since $\partial e^{\intQ} = Q_{(-1)} e^{\intQ}$, we have
  \[
   \sum_{j=0}^{k} Q_{(-j-1)} e^{\intQ}_{(-k+j)} \One = k e^{\intQ}_{(-k-1)} \One
  \]
  for $k \ge 1$.  Therefore, we obtain
  
  \begin{multline*}
  e^{\intQ}_{(0)} \Bigl( \tilR{1} \tilR{0}^m \One - \frac{(m+1)(\K-1) + 1}{(m+1)(\K-1)} \tilR{0}^{m+1} \One \Bigr) \\
   = \sum_{i=2}^{m+1} (-1)^i (\K-1)^{i-1} \frac{m!}{(m-i+1)!} (i-1) \tilR{0}^{m-i+1} e^{\intQ}_{(-i)} \One \\
   + \sum_{k=1}^{m} \sum_{j=0}^{k} 
   (-1)^{k} (\K-1)^{k} \frac{m!}{(m-k)!} \tilR{0}^{m-k} Q_{(-j-1)} e^{\intQ}_{(-k+j)} \One \\
   = \sum_{i=2}^{m+1} (-1)^i (\K-1)^{i-1} \frac{m!}{(m-i+1)!} (i-1) \tilR{0}^{m-i+1} e^{\intQ}_{(-i)} \One \\
   + \sum_{k=1}^{m} (-1)^k (\K-1)^k \frac{m!}{(m-k)!} \tilR{0}^{m-k} k e^{\intQ}_{(-k-1)} \One = 0.
  \end{multline*}
   \end{proof}
\end{lemma}

\begin{proposition}
 \label{prop:2}
 We have $e^{\intQ}_{(0)} \gfW(u) = 0$.
\begin{proof}
 We split terms in the definition of $\gfW(u)$ into two parts; terms with the factor $\tilR{1}$ and terms without
 $\tilR{1}$. Note that the screening operator $e^{\intQ}_{(0)}$ commutes with $\tilR{i}$ for all $i \ne 0$, $1$.
 Then we have
 \begin{multline*}
  e^{\intQ}_{(0)} \gfW(u) = e^{\intQ}_{(0)} \sum_{m=0}^{N} (-1)^{m} \prod_{j=1}^{m} \frac{j(\K-1) + 1}{j(\K-1)} 
  \hspace{-0.3cm} \sum_{N \ge i_1 > \dots > i_{N-m} \ge 1} \hspace{-0.7cm} 
  \tilR{i_1} \cdots \tilR{i_{N-m}} \tilR{0}^{m} \One \\
  = e^{\intQ}_{(0)} \Bigl\{
  \sum_{m=0}^{N} (-1)^{m} \prod_{j=1}^{m} \frac{j(\K-1) + 1}{j(\K-1)} 
  \hspace{-0.3cm} \sum_{N \ge i_1 > \dots > i_{N-m-1} \ge 2} \hspace{-0.7cm} 
  \tilR{i_1} \cdots \tilR{i_{N-m-1}} \tilR{1} \tilR{0}^{m} \One \\
+ \sum_{m=0}^{N} (-1)^{m+1} \prod_{j=1}^{m+1} \frac{j(\K-1) + 1}{j(\K-1)} 
  \hspace{-0.3cm} \sum_{N \ge i_1 > \dots > i_{N-m-1} \ge 2} \hspace{-0.7cm} 
  \tilR{i_1} \cdots \tilR{i_{N-m-1}} \tilR{0}^{m+1} \One \Bigr\} = 0
 \end{multline*}
 by \reflemma{lemma:1}.
\end{proof}
\end{proposition}

By \refprop{prop:1} and \refprop{prop:2}, we have 
$\W'_m \in \Ker e^{\intQ}_{(0)} \cap \bigcap_{i=1}^{N-1} \Ker e^{\intA_i}_{(0)}$ for $m=1$, $\dots$, $N$.
Thus, we have the following proposition by a consequence of \refprop{prop:G-isom}

\begin{proposition}
 \label{prop:W-belong}
 For $m=1$, $\dots$, $N$, we have $\W'_m \in \WalgR$.
\end{proposition}

To construct elements which give strong generators of the vertex algebra $\WalgR$, we need
to normalize the elements $\W'_m$ for $m=2$, $\dots$, $N$. 

\begin{lemma}
 \label{lemma:2}
 Let $q$ be an indeterminate. For $l$, $m \ge k$, we have the following identity:
 \[
 \sum_{l=k}^{N} (-1)^{l-k} \binom{N-m+k}{l} \binom{l}{k} \prod_{j=k+1}^{l} \Bigl( 1 + \frac{q}{j} \Bigr)
 = \frac{(-1)^{N-m}}{(N-m)!} \prod_{j=0}^{N-m-1} (q - j).
 \]
 \begin{proof}
  By direct calculation, we have
\begin{multline*}
\sum_{l=k}^{N} (-1)^{l-k} \binom{N-m+k}{l} \binom{l}{k} \prod_{j=k+1}^{l} \Bigl( 1 + \frac{q}{j} \Bigr) \\
 = \sum_{l=0}^{N-m} (-1)^l \binom{N-m+k}{l+k} \frac{1}{l!} \prod_{j=1}^{l} (q+k+j) \\
 = \frac{(N-m+k)!}{(N-m)!k!} \sum_{l=0}^{N-m} \frac{(-N+m)_{l} (q+k+1)_{l}}{(k+1)_{l}} \frac{1^l}{l!}
\end{multline*}
where $(x)_{l} = \prod_{j=0}^{l-1} x+j$. The RHS can be described by the hypergeometric function
${}_2 F_1(a, b, c; z)$, and thus by the $\Gamma$-function $\Gamma(z)$ with applying Gauss's 
hypergeometric theorem. Then we have
\begin{multline*}
 \sum_{l=k}^{N} (-1)^{l-k} \binom{N-m+k}{l} \binom{l}{k} \prod_{j=k+1}^{l} \Bigl( 1 + \frac{q}{j} \Bigr) \\
 = \frac{(N-m+k)!}{(N-m)!k!} {}_2 F_1(-N+m, q+k+1, k+1; 1) \\
 = \frac{(N-m+k)!}{(N-m)!k!} \frac{\Gamma(k+1) \Gamma(N-m-q)}{\Gamma(N-m+k+1)\Gamma(-q)} \\
 = \frac{(N-m+k)!}{(N-m)!k!} \frac{k!}{(N-m+k)!} \prod_{j=0}^{N-m-1} (-q+j) 
 = \frac{(-1)^{N-m}}{(N-m)!} \prod_{j=0}^{N-m-1} (q - j).
\end{multline*}
 \end{proof}
\end{lemma}

\begin{lemma}
 \label{lemma:3}
 Let $\OpX_0$, $\OpX_1$, $\dots$, $\OpX_N$ be operators. For $m=1$, $\dots$, $N$, we have
\begin{multline*}
 \sum_{k=0}^{N} (-1)^k \prod_{j=1}^{k} \frac{j(\K-1)+1}{j(\K-1)} 
 \sum_{1 \le i_1 < \dots < i_{N-k} \le N} \hspace*{-5mm} 
 e_m(\overbrace{\OpX_0, \dots, \OpX_0}^{k\text{-times}}, \OpX_{i_1}, \dots, \OpX_{i_{N-k}})
 \\
 = \prod_{j=1}^{N-m} \frac{j(\K-1) - \K}{j(\K-1)} \sum_{k=0}^{m} (-1)^k 
 \Bigl(\prod_{j=1}^{k} \frac{j(\K-1) + 1}{j(\K-1)} \Bigr)
 e_{m-k}(\OpX_1, \dots, \OpX_N) \OpX_0^k.
\end{multline*}
\begin{proof}
Below we write $e_{m-k}(\OpX_1, \dots, \OpX_N)$ by $e_{m-k}$ for short.
\begin{multline*}
 \sum_{k=0}^{N} (-1)^k \prod_{j=1}^{k} \frac{j(\K-1)+1}{j(\K-1)} 
 \sum_{1 \le i_1 < \dots < i_{N-k} \le N} \hspace*{-5mm} e_m(\overbrace{\OpX_0, \dots, \OpX_0}^{k\text{-times}}, \OpX_{i_1}, \dots, \OpX_{i_{N-k}})
 \\
 = \sum_{l=0}^{N} (-1)^l \prod_{j=1}^{l} \frac{j(\K-1)+1}{j(\K-1)} \hspace*{-3mm} \sum_{i_1 > \dots > i_{N-l}}
 \hspace*{-2mm}
 \sum_{k=0}^{\min(l, m)} \hspace*{-2mm} \sum_{j_1 < \dots < j_{m-k}} \hspace*{-4mm} \OpX_{i_{j_1}} \cdots \OpX_{i_{j_{m-k}}} 
 \binom{l}{k} \OpX_0^k \\
= \sum_{l=0}^{N} (-1)^l \prod_{j=1}^{l} \frac{j(\K-1)+1}{j(\K-1)} 
 \sum_{k=0}^{\min(l, m)} \hspace*{-2mm} \sum_{i_1 > \dots > i_{m-k}} \hspace*{-4mm} 
 \OpX_{i_{1}} \cdots \OpX_{i_{m-k}} \binom{l}{k} \OpX_0^k \\
\shoveleft{= \sum_{k=0}^{m} (-1)^k \Bigl( \prod_{j=1}^{k} \frac{j(\K-1)+1}{j(\K-1)} \Bigr)} \\
 \cdot \sum_{l=k}^{N} (-1)^{l-k} \binom{N-m+k}{l} \binom{l}{k} \prod_{j=k+1}^{l} \Bigl(1 + \frac{1}{j(\K-1)} \Bigr)
 e_{m-k} \OpX_0^k
\end{multline*}
 Applying \reflemma{lemma:2} for $q = 1/(\K-1)$, we obtain the identity of the lemma.
\end{proof}
\end{lemma}

Applying the above lemma for $\OpX_i = \F{i}$, we have identities for 
the elements $\W'_m$ ($m=1$, $\dots$, $N$).

\begin{multline*}
 \W'_m = \sum_{k=0}^{N} (-1)^k \Bigl( \prod_{j=1}^{k} \frac{j(\K-1)+1}{j(\K-1)} \Bigr)
 \sum_{i_1 < \dots < i_{N-k}} e_m(\overbrace{\F{0}, \dots, \F{0}}^{k\text{-times}}, \F{i_1}, \dots, \F{i_{N-k}}) \One \\
 = \prod_{j=1}^{N-m} \frac{j(\K-1)-\K}{j(\K-1)} \sum_{k=0}^{m} (-1)^k 
 \Bigl(\prod_{j=1}^{k} \frac{j(\K-1)+1}{j(\K-1)}\Bigr) 
 e_{m-k}(\F{1}, \dots, \F{N}) \F{0}^k \One
\end{multline*}

\begin{definition}
 \label{def:generators}
 For $m=1$, $\dots$, $N$, set
 \begin{multline}
  \label{eq:5}
  \W''_m = (-1)^m \prod_{j=1}^{N-m} \frac{j(\K-1)}{j(\K-1)-\K} \W'_m \\
  = \sum_{k=0}^{m} (-1)^{m+k} 
 \Bigl(\prod_{j=1}^{k} \frac{j(\K-1)+1}{j(\K-1)}\Bigr) 
\cdot e_{m-k}(\F{1}, \dots, \F{N}) \F{0}^k \One \in \WalgR.
 \end{multline}
 Then, we define elements $\W_m \in \WalgR$ for $m=1$, $\dots$, $N$ inductively as follows:
\begin{equation}
 \label{eq:6} 
 \W_m = \W''_m - \sum_{k=0}^{m-1} (-1)^k \Bigl(\prod_{j=1}^{k} \frac{j(\K-1)+1}{j(\K-1)}\Bigr) \cdot \W_{m-k (-1)}
  (H_{(-1)})^k \One,
\end{equation}
 and $\W_1 = 0$. 
\end{definition}

In \cite[Lemma 2.3.5]{FS}, a conformal vector $\omega$ of the vertex algebra $\WalgCK \otimes \C[\K, \K^{-1}]$ is
explicitly given over $\C[\K, \K^{-1}]$. By direct calculation, we have the following 
relations between $\W_1$, $\W_2$ and $H$, $\omega$.

\begin{proposition}
 \label{prop:conf-vec}
 We have the following identities between the elements $\W_2$ and $H$, $\omega$;
\begin{equation*}
 \omega = - \frac{1}{\K} \W_2 
 + \frac{N-1}{2} \frac{2\K-1}{(\K-1)^2} H_{(-1)} H + \frac{N}{2} \frac{(2N-3)\K - (2N-2)}{\K-1} \partial H,
\end{equation*}
 over $\C[\K, \K^{-1}, (\K-1)^{-1}]$.
\end{proposition}

In the rest of this section, we show that the $N+1$ elements $E$, $F$, $H$ and $\W_m$ for $m = 2$, 
$\dots$, $N-1$ strongly generate the subregular $\calW$-algebra $\WalgR \otimes \C_{\Kv}$ for 
$\Kv \in \C \setminus \{1\}$.

For a vertex algebra $V$, let $\Ctwo{V} = V / C_2(V)$ be Zhu's $C_2$ Poisson algebra of $V$, where
$C_2(V) = V_{(-2)} V$. 
For an element $a \in V$, we denote its image in $\Ctwo{V}$ by $\overline{a} \in \Ctwo{V}$.

\begin{lemma}
 \label{lemma:lead-term}
 For arbitrary $\Kv \in \C \setminus \{1\}$, 
 we have $\barW_m = e_m(\barmu_1, \dots, \barmu_N)$ in $\Ctwo{\calV^{\Kv}}$ for $m=1$, $\dots$, $N$.
\begin{proof}
 Since $Y = e^{\intY}_{(-2)} e^{- \intY} \equiv 0$ modulo $C_2(\calV^{\Kv})$, we have $\barH = \barmu_0$.
 Thus, by \refeq{eq:5}, we have
 \begin{equation}
\label{eq:barW-C2} 
 \barW''_m 
  = \sum_{k=0}^{m} (-1)^k
 \Bigl(\prod_{j=1}^{k} \frac{j(\K-1)+1}{j(\K-1)}\Bigr) \cdot e_{m-k}(\barmu_1, \dots, \barmu_N) \,\barH^k.
\end{equation}
 Then, the claim of the lemma follows from \refeq{eq:barW-C2} and \refeq{eq:6} by induction on $m$.
\end{proof}
\end{lemma}

\begin{lemma}
 \label{lemma:alg-indep}
 For arbitrary $\Kv \in \C \setminus \{1\}$, 
 the elements $\barH$, $\barW_2$, $\dots$, $\barW_N \in \Ctwo{\calW^{\Kv}}$
 are algebraically independent over the field $\C$.
\begin{proof}
 First, note that it is enough to see that $\barH$, $\barW_2$, $\dots$, $\barW_N$ are
 algebraically independent in $\Ctwo{\calV^{\Kv}}$, since we have $C_2(\calW^{\Kv}) \subset C_2(\calV^{\Kv})$.
 Below we write $\bare_m = e_m(\barmu_1, \dots, \barmu_N)$ for short.
 Note that we only have elements of form $e^{m \intY}$ for $m \in \Z$ in the lattice part of $\ValgR$,
 and thus $\overline{Q}$, $\overline{A}_1$, $\dots$, $\overline{A}_{N-1}$ are linearly independent in
 $\Ctwo{\calV^{\Kv}}$, while $Y = e^{\intY}_{(-2)} e^{- \intY} \equiv 0$ modulo 
 $C_2(\calV^{\Kv})$.
 Since $\barmu_i \in \bigoplus_{j=1}^{N-1} \C \overline{A}_j$ for all $i=1$, $\dots$, $N$ and 
 $\barH = \barmu_0 \not\in \bigoplus_{j=1}^{N-1} \C \overline{A}_j$, $\barH$ is algebraically
 independent of $\barW_2 = \bare_2$, $\dots$, $\barW_N = \bare_N$.

 We identify the vector space $\bigoplus_{j=1}^{N-1} \C \overline{A}_j$ with the Cartan subalgebra
\[
 \frh = \Bigl\{\, \sum_{i=1}^{N} c_i \varepsilon_i \in \bigoplus_{i=1}^{N} \C \varepsilon_i \,\Bigm|\, 
 c_1 + \dots + c_N = 0 \, \Bigr\} \simeq \C^{N-1}
\]
 by the standard way; $\overline{A}_j = \varepsilon_j - \varepsilon_{j+1}$ ($j=1$, $\dots$, $N-1$). 
 Under this identification, we have $\barmu_i = \varepsilon_i - (1/N) \sum_{j=1}^{N} \varepsilon_j$ for 
 $i=1$, $\dots$, $N$. Then, the symmetric polynomials $\bare_2$, $\dots$, $\bare_N$ are algebraically
 independent and we have $\C[\frh]^{\frS_N} = \C[\bare_2, \dots, \bare_N]$, while $\bare_1 = 0$
 by the classical fact on the Weyl-group-invariant subalgebra $\C[\frh]^{\frS_N}$. 
 Thus, we have the assertion of the lemma.
\end{proof}
\end{lemma}

\begin{proposition}
 \label{prop:Poisson-center}
 The Poisson center of $C_2$ Poisson algebra $\Ctwo{\calW^{\Kv}}$ is generated by $\barW_2$, $\dots$, $\barW_N$.
 \begin{proof}
  Since $(\Mu_i, Y) = 0$ for $i=1$, $\dots$, $N$, it is easy to check that 
  $\{\barmu_i, \barE\} = \{\barmu_i, \barF\} = 0$. Thus, $\bare_2$, $\dots$, $\bare_N$ are Poisson
  central in $\Ctwo{\calV^{\Kv}}$. By \cite[Lemma 6.12]{Genra2}, 
  $\FMS \circ \Wak \circ \mu: \calW^{\Kv} \longrightarrow \calV^{\Kv}$ induces an embedding
  $\Ctwo{\calW^{\Kv}} \hookrightarrow \Ctwo{\calV^{\Kv}}$, and thus $\bare_2$, $\dots$, $\bare_N$
  are Poisson central also in $\Ctwo{\calW^{\Kv}}$. By a consequence of \cite{DK}, the Poisson
  center of $\Ctwo{\calW^{\Kv}}$ is isomorphic to $\C[\frh]^{W}$, and $\barmu_2$, $\dots$, $\barmu_N$
  are algebraically independent by \reflemma{lemma:alg-indep}. Therefore, the Poisson center
  of $\Ctwo{\calW^{\Kv}}$ is $\C[\barW_2, \dots, \barW_N] \simeq \C[\frh]^{W}$.
 \end{proof}
\end{proposition}

It is easy to check that the elements $E$, $H$, $\W_2$, $\W_3$, $\dots$, $\W_{N-1}$, $F$ have conformal weights
$1$, $1$, $2$, $3$, $\dots$, $N-1$, $N-1$ respectively. Note that, for arbitrary $\Kv \in \C$, 
the vertex algebra $\WalgCK \otimes \C_{\Kv}$ is of type $\calW(1, 1, 2, 3, \dots, N-1, N-1)$,
i.e. $\WalgCK \otimes \C_{\Kv}$ has $N+1$ strong generators of conformal weight $1$, $1$, $2$, $3$, $\dots$, $N-1$, $N-1$.

\begin{theorem}
 \label{thm:strong-gen}
 For arbitrary $\Kv \in \C \setminus \{1\}$, 
 the elements $E$, $H$, $\W_2$, $\dots$, $\W_{N-1}$, $F$ strongly generate the vertex algebra 
 $\WalgR \otimes \C_{\Kv} = \calW^{\Kv}$.
\begin{proof}
 The vertex algebra $V = \WalgR \otimes \C_{\Kv}$ is decomposed as $V = \bigoplus_{d \ge 0} V_d$ where
 $V_d$ is the subspace of conformal weight $d$. 
 Since $V$ is of type $\calW(1,1,2, \dots, N-1, N-1)$, $V_1$ is two-dimensional. On the other hand, 
 there exist two linearly independent elements $E$ and $H$, and thus we have $V_1 = \C E \oplus \C H$.

 By induction on $d$, we show that 
\begin{equation}
\label{eq:3} 
  V_d \subset \Span \{ \EleX_{1 (-n_1)} \cdots \EleX_{k (-n_k)} \One \,|\, \EleX_i = E, H, \W_2, \dots, \W_d, n_i \ge 1 \}
\end{equation}
 for $d = 2$, $\dots$, $N-1$. Assume that \refeq{eq:3} holds for $V_{d'}$ with $d' \le d-1$. Set
 \[
  U_d = V_d \cap \Span \{ \EleX_{1 (-n_1)} \cdots \EleX_{k (-n_k)} \One \,|\, \EleX_i = E, H, \W_2, \dots, \W_{d-1}, n_i \ge 1 \}.
 \]
 Then, $U_d$ is codimension one in $V_d$ since $V$ has exactly one strong generator of conformal weight
 $d$. We show that $\W_d \not\in U_d$. Indeed, assume that we have an identity
\begin{equation}
\label{eq:4} 
 \W_d = \sum_{p} c^{(p)} \EleX^{(p)}_{1 (-n_1)} \cdots \EleX^{(p)}_{k_p (-n_{k_p})} \One
\end{equation}
 where $c^{(p)} \in \C$,  $\EleX^{(p)}_{i} = E$, $H$, $\W_2$, $\dots$, $\W_{d-1}$ and $n_i \ge 1$. 
 Terms containing $E_{(-n)}$ ($n \ge 1$) have positive $H_{(0)}$-eigenvalues while $H_{(0)} \W_m = 0$
 for all $m$ and $H_{(0)} H = 0$. By using decomposition into $H_{(0)}$-eigenspaces, we may assume that
 the identity \refeq{eq:4} holds for $\EleX_i^{(p)} = H$, $\W_2$, $\dots$, $\W_{d-1}$. 
 Taking modulo $C_2(V)$, we have an identity 
 $\barW_d = \sum_{p} c^{(p)} \overline{\EleX}^{(p)}_{1 (-n_1)} \cdots \overline{\EleX}^{(p)}_{k_p (-n_{k_p})} \One$ in
 the algebra $\Ctwo{V}$. It contradicts \reflemma{lemma:alg-indep}, and hence we have
 $\W_d \not\in U_d$. Therefore, we have $V_d = U_d \oplus \C \W_d$ and the induction completes.

Similarly to the above, setting
 \[
  U_{N-1} = V_{N-1} \cap \Span \{ \EleX_{1 (-n_1)} \cdots \EleX_{k (-n_k)} \One \,|\, \EleX_i = E, H, \W_2, \dots, \W_{N-1}, n_i \ge 1 \},
 \]
 we have $\Dim V_{N-1} / U_{N-1} = 2$. Since the element $\barW_{N-1}$ is algebraically independent of 
 $\barH$, $\barW_2$, $\dots$, $\barW_{N-2}$ in $\Ctwo{V}$ by \reflemma{lemma:alg-indep}, we have
 $\W_{N-1} \not\in U_{N-1}$. Note that the element $F$ has weight $-1$ with respect to the action
 of $H_{(0)}$, and thus we have $F \not\in U_{N-1} \oplus \C \W_{N-1}$. Therefore we have
 $V_{N-1} = U_{N-1} \oplus \C \W_{N-1} \oplus \C F$.

 Since $V$ have no strong generator with conformal weight bigger than $N-1$, we have
 \[
  V = \Span \{ \EleX_{1 (-n_1)} \cdots \EleX_{k (-n_k)} \One \,|\, \EleX_i = E, H, \W_2, \dots, \W_{N-1}, F, n_i \ge 1 \}.
 \]
 It implies that $E$, $H$, $\W_2$, $\dots$, $\W_{N-1}$, $F$ strongly generates $V$.
\end{proof}
\end{theorem}

\begin{remark}
 Consider another set of elements 
 \[
  U_m  = \sum_{k=0}^{m} (-1)^{m+k} 
 \Bigl(\prod_{j=1}^{k} \frac{j(\K-1)+1}{j(\K-1)}\Bigr) 
\cdot e_{m-k}(\R{1}, \dots, \R{N}) \R{0}^k \One
 \]
 of the vertex algebra $\ValgCK$ instead of $\W_m$ for $m=1$, $\dots$, $N$, 
 where $\R{i} = \F{i} + H_{(-1)}$ for $i=0$, $1$, $\dots$, $N$.
 Since $\R{0} = (\K-1) (\partial + Y_{(-1)})$, $U_m$ is well-defined even over the ring $\C[\K]$
 for $m=1$, $\dots$, $N$. One can easily check that we have 
 $e^{\intQ}_{(0)} U_m = e^{\intA_i}_{(0)} U_m = 0$ for all $i=1$, $\dots$, $N-1$ and $m=1$, $\dots$,
 $N$ merely by replacing $\F{i}$ by $\R{i}$ for $i=0$, $\dots$, $N$ 
 in the proofs of lemmas and propositions in this section, since $H$ lies in the kernels of the
 screening operators. 
 Thus, $U_1$, $\dots$, $U_N$ are elements of $\WalgCK$.

 Moreover, the proofs of \reflemma{lemma:alg-indep} and \refthm{thm:strong-gen} work also analogously
 for $U_2$, $\dots$, $U_{N-1}$. Therefore, we have another set of strong generators
 $E$, $H$, $U_2$, $\dots$, $U_{N-1}$, $F$ of the vertex algebra $\WalgCK \otimes \C_{\Kv} = \calW^{\Kv}$ for
 arbitrary $\Kv \in \C$. While this set of elements gives strong generators even for $\Kv = 1$,
 we prefer the elements $\W_m$ ($m=2$, $\dots$, $N$) because $\W_2$, $\dots$, $\W_N$ strongly 
 generate the center of the vertex algebra $\WalgCK \otimes \C_{0} = \calW^{0} = \calW^{-N}(\frsl_N, f_{sub})$
 at critical level as we discuss in the following section.
\end{remark}

\section{Structure of the subregular $\calW$-algebra at critical level}
\label{sec:critical-level}

In this section, we consider the strong generators $E$, $H$, $\W_2$, $\dots$, $\W_{N-1}$, $F$
of the subregular $\calW$-algebra at critical level ($\K = 0$), and study the OPEs between 
these generators.
Throughout this section, we specialize $\K$ to $0$, and consider the elements $E$, $F$, $H$, 
$\W_2$, $\dots$, $\W_{N}$ in the vertex algebra $\calW^{-N}(\frsl_N, f_{sub}) = \WalgR \otimes \C_0$.

First, by \refdef{def:generators} and \refeq{eq:def-Wm1}, we have
\[
 \W_m = (-1)^m \W'_m = e_m(\partial + \Mu_{1 (-1)}, \dots, \partial + \Mu_{N (-1)}) \One
\]
for $m=2$, $\dots$, $N$. 
Since $\Mu_i \in \bigoplus_{j=1}^{N-1} \C A_{j}$ and $(A_j, \blkbar) = 0$ for $i=1$, $\dots$, $N$
and $j=1$, $\dots$, $N-1$, the element $\W_m$ is central for $m=2$, $\dots$, $N$.

\begin{proposition}
 \label{prop:gen-center}
 The elements $\W_2$, $\dots$, $\W_N$ strongly generate the center of the vertex algebra
 $\calW^{-N}(\frsl_N, f_{sub})$.
\begin{proof}
 By \cite[Theorem 1.1]{Arakawa}, the center of the vertex algebra $\calW^{-N}(\frsl_N, f_{sub})$ coincides
 with the center of the universal affine vertex algebra $V^{-N}(\frsl_N)$, and it is $\Z_{\ge 0}$-graded.
 Hence, it is the vertex algebra of type $\calW(2, 3, \dots, N)$. Note that the elements $\barW_2$, $\dots$,
 $\barW_N$ are algebraically independent in the $C_2$ Poisson algebra by \reflemma{lemma:alg-indep}.
 Thus, applying the same argument in the proof of \refthm{thm:strong-gen} to the elements $\W_2$, 
 $\dots$, $\W_N$ of the center, the assertion of the proposition holds.
\end{proof}
\end{proposition}

Note that we know the OPEs
\begin{equation}
\label{eq:OPE-H} 
 H(z) E(w) \sim \frac{1}{z-w} E(w), \qquad H(z) F(w) \sim \frac{-1}{z-w} F(w).
\end{equation}
To describe complete structure of the vertex algebra $\calW^{-N}(\frsl_N, f_{sub})$ algebraically,
we discuss the OPE between $E$ and $F$. First, by direct computation, we have the following lemma.

\begin{lemma}
 \label{lemma:5}
 We have the commutation relation $[e^{\XI}_{(m)}, \R{i}] = (- m - 1) e^{\XI}_{(m-1)}$
 for $i=1$, $\dots$, $N$, where $\R{i} = - \partial + H_{(-1)} - \Mu_{i (-1)}$ is the
 operator defined in \refsec{sec:strong-gen}.
\end{lemma}

In the following, we set $\W_0 = \One$ and $\W_1 = 0$. 
The following lemma is essentially due to Molev. See \cite[Proposition 12.4.4]{Molev}.

\begin{lemma}
 \label{lemma:Molev}
 For $m = 1$, $\dots$, $N$, we have the following identity
\begin{multline*}
 e_{m}(\partial - H_{(-1)} + \Mu_{1 (-1)}, \dots, \partial - H_{(-1)} + \Mu_{N (-1)}) \\
= \sum_{k=0}^{m} \binom{N-k}{m-k} \W_{k (-1)} (\partial - H_{(-1)})^{m-k}.
\end{multline*}
\begin{proof}
 Define $\zeta_m \in \C[\Mu_{i (-n)} \,|\, i=1, \dots, N, n \ge 1]$ by
\[
 \sum_{m=0}^{N} \zeta_m \partial^m = (\partial + \Mu_{N (-1)}) \cdots (\partial + \Mu_{1 (-1)}).
\]
Let $u$ be an indeterminate, and we obtain the following by replacing $\partial$ by $u + \partial$,
\begin{multline*}
 (u + \partial + \Mu_{N (-1)}) \cdots (u + \partial + \Mu_{1 (-1)})
 = \sum_{k=0}^{N} \zeta_k (u + \partial)^{N-k} \\
 = \sum_{m=0}^{N} \sum_{k=0}^{m} \binom{N-k}{m-k} \zeta_k \partial^{m-k} u^{N-m}.
\end{multline*}
On the other hand, the LHS is clearly equal to 
\[
 \sum_{m=0}^{N} e_m(\partial + \Mu_{1 (-1)}, \dots, \partial+\Mu_{N(-1)}) u^{N-m},
\]
and thus we have $\zeta_m \One = e_m(\partial + \Mu_{1 (-1)}, \dots, \partial+\Mu_{N(-1)}) \One = \W_m$.
Since $\Mu_1$, $\dots$, $\Mu_N$ are central, we have $\zeta_m = \W_{m (-1)}$.

Since $H_{(-1)}$ commutes with $\Mu_{i (-1)}$ for $i=1$, $\dots$, $N$, one can replace 
$\partial$ by $\partial - H_{(-1)}$ in the above identity, thus we have
\begin{multline*}
 e_m(\partial - H_{(-1)} + \Mu_{1 (-1)}, \dots, \partial - H_{(-1)} + \Mu_{N(-1)}) \\
 = \sum_{m=0}^{N} \sum_{k=0}^{m} \binom{N-k}{m-k} \W_{k (-1)} (\partial - H_{(-1)})^{m-k}
\end{multline*}
 as the coefficient of $u^{N-m}$ for $m=1$, $\dots$, $N$. 
\end{proof}
\end{lemma}

Now we describe the OPE between $E$ and $F$ in terms of our strong generators.

\begin{theorem}
 \label{thm:OPE}
 We have the following OPE:
 \[
  E(z) F(w) \sim (-1)^{N+1} \sum_{n=1}^{N-1} \frac{n!}{(z-w)^n} \sum_{m=0}^{N-n} \binom{N-m}{n}
 \Bigl( \W_{m (-1)} (\partial - H_{(-1)})^{N-n-m} \One \Bigr)(w)
 \]
\begin{proof}
 For $m \le N+1$, we have by using \reflemma{lemma:5} repeatedly
\begin{multline*}
 E_{(N-m)} F = - e^{\XI}_{(N-m)} \R{N} \cdots \R{1} e^{- \XI} \\
 = - \R{N} \cdots \R{1} e^{\XI}_{(N-m)} e^{-\XI} 
 - \dots
 - \Bigl(\prod_{j=1}^{N-m} (-j-1) \Bigr) \sum_{i_1 > \dots > i_{m}} \R{i_1} \cdots \R{i_m} e^{\XI}_{(0)} e^{- \XI} \\
 + \Bigl(\prod_{j=1}^{N-m} (-j-1) \Bigr) \sum_{i_1 > \dots > i_{m-1}} \R{i_1} \cdots \R{i_{m-1}} e^{\XI}_{(-1)} e^{- \XI} \\
 = (-1)^{N+1} (N-m+1)! e_{m-1}(\partial - H_{(-1)} + \Mu_{1 (-1)}, \dots, \partial - H_{(-1)} + \Mu_{N (-1)}) \One.
\end{multline*}
Here we used $e^{\XI}_{(n)} e^{- \XI} = 0$ for $n \ge 0$ and $e^{\XI}_{(-1)}$ commutes with $\R{i}$ ($i=1$, $\dots$, 
 $N$). Applying \reflemma{lemma:Molev}, we obtain
\begin{equation}
\label{eq:n-prod-E-F} 
 E_{(N-m)} F = (-1)^{N+1} (N-m+1)! \sum_{k=0}^{m-1} \binom{N-k}{m-k-1} \W_{k (-1)} (\partial - H_{(-1)})^{m-k-1} \One,
\end{equation}
and the OPE of the theorem is an immediate consequence of it.
\end{proof}
\end{theorem}

Using the strong generators $E$, $F$, $H$ and $\W_m$ ($m=2$, $\dots$, $N-1$), we can determine
the structure of the Zhu algebra of $\calW^{-N}(\frsl_N, f_{sub})$ explicitly.

Let $V = \bigoplus_{\Delta \ge 0} V_{\Delta}$ be a $\Z_{\ge 0}$-graded vertex algebra, and
we denote the degree of a homogeneous element $\EleX \in V$ by $\Delta(\EleX)$. For 
For $a \in V_{\Delta}$ and $b \in V$, we define
\[
 a \circ b = \sum_{j=0}^{\Delta} \binom{\Delta}{j} a_{(j-2)} b, \qquad
 a * b = \sum_{j=0}^{\Delta} \binom{\Delta}{j} a_{(j-1)} b.
\] 
Then, the vector space $\Zhu{V} = V / (V \circ V)$ has a structure of an associative algebra
by the multiplication induced by $*$, called the Zhu algebra of $V$ \cite{Zhu, FZ, DK}.
For $V = \calW^{\Kv-N}(\frsl_N, f_{sub})$, the Zhu algebra 
$\Zhu{\calW^{\Kv-N}(\frsl_N, f_{sub})} = U(\frsl_N, f_{sub})$ is known as the finite 
$\calW$-algebra associated with $\frsl_N$ and $f_{sub}$ by the result of De Sole and Kac in
\cite{DK}, and in particular it does not depend on the level $\Kv - N$. See also \cite{Arakawa2}. Moreover, in \cite{Premet},
Premet showed that the finite $\calW$-algebra $U(\frsl_N, f_{sub})$ was isomorphic to Smith's algebra
introduced by Smith in \cite{Smith}. Below we describe the structure of $\Zhu{\calW^{-N}(\frsl_N, f_{sub})}$.

The vertex algebra $\calV^0$ is $\Z$-graded by $\Delta(A_i) = 1$ for $i=1$, $\dots$, $N-1$,
$\Delta(Q) = 1$, $\Delta(Y) = 1$ and $\Delta(e^{\pm \intY}) = \pm 1$. This grading induces
$\Z_{\ge 0}$-grading on $\calW^{-N}(\frsl_N, f_{sub}) = \calW^{0} \subset \calV^0$; i.e. we
have $\Delta(E) = 1$, $\Delta(F) = N-1$, $\Delta(H) = 1$ and $\Delta(\W_m) = m$ for $m=2$, $\dots$,
$N$. Note that the grading coincides with the conformal weight in spite of $\calW^{-N}(\frsl_N, f_{sub})$
is not a vertex operator algebra.

First, we discuss the $C_2$ Poisson algebra $\overline{A} := \Ctwo{\calW^{-N}(\frsl_N, f_{sub})}$. The elements
$\barW_2$, $\dots$, $\barW_N$ are Poisson central and algebraically independent by \refprop{prop:gen-center}
and \reflemma{lemma:alg-indep}, and thus $\overline{A}$ is
a Poisson algebra over $\C[\barW_2, \dots, \barW_{N-1}]$. By \refeq{eq:n-prod-E-F} for $m=N+1$, we have 
\[
 \barE \barF = \sum_{k=0}^{N} (-1)^{k+1} \barW_{k} \barH^{N-k}.
\]
By the result of \cite{DK} or \cite{Arakawa2}, the Poisson algebra $\overline{A}$ is the graded algebra of $U(\frsl_N, f_{sub})$, and thus is the coordinate ring of the Slodowy
slice $\bbS$ in $\frsl_N$ associated with $f_{sub}$, which is known as the simultaneous 
deformation of the Kleinian singularity of type $A_{N-1}$. These fact gives an isomorphism of
commutative algebras
\[
 \overline{A} = \C[\barE, \barF, \barH, \barW_m \,|\, m=2, \dots, N-1] \bigm/ 
\bigl(\barE \barF - \sum_{k=0}^{N} (-1)^{k+1} \barW_{k} \barH^{N-k} \bigr),
\]
and the subalgebra $\C[\W_2, \dots, \W_N]$ is the Poisson center of $\overline{A}$.
Note that the Poisson brackets between these elements are given by $\{\barH, \barE\} = \barE$,
$\{\barH, \barF\} = - \barF$, and 
$\{\barE, \barF\} = \sum_{k=0}^{N} (-1)^{k+1} (k+1) \barW_{k} \barH^{N-k-1}$ by \refeq{eq:n-prod-E-F}.

In \cite{Smith}, Smith studied an associative algebra $R$ generated by the elements $x$, $y$ and $h$ 
subject to the relation $[h, x] = x$, $[h, y] = - y$, $[x, y] = f(h)$ where $f$ is a polynomial of
degree $N-1$.
The Zhu algebra $A := \Zhu{\calW^{-N}(\frsl_N, f_{sub})}$ is generated by the image of $E$, $F$, $H$
and $\W_m$ ($m=2$, $\dots$, $N-1$) by \refthm{thm:strong-gen}. Note that, by \refprop{prop:gen-center},
the center of $A$ is $\C[\W_2, \dots, \W_{N}]$. As an associative algebra
over $\C[\W_2, \dots, \W_{N-1}]$, $A$ is generated by $E$, $F$ and $H$. 
It is easy to see that $[H, E] = E$ and $[H, F] = - F$ by \refeq{eq:OPE-H}. Using \refeq{eq:n-prod-E-F}
and the skew-symmetry $F_{(n)} E = \sum_{j \ge 0} (-1)^{n+j-1} \partial^j (E_{(n+j)} F)/j!$,
$[E, F]$ is a polynomial in $H$ of degree $N-1$ with coefficients in $\C[\W_2, \dots, \W_{N-1}]$.
Thus, the Zhu algebra $A$ is Smith's algebra over 
$\C[\W_2, \dots, \W_{N-1}]$. The isomorphism between the finite $\calW$-algebra 
$U(\frsl_N, f_{sub})$ and Smith's algebra is a well-known result by Premet \cite[Theorem 7.10]{Premet}.

The vertex algebra $\calW^{-N}(\frsl_N, f_{sub})$ can be realized also by using fermionic fields.
The realization was conjectured by Adamovi\'{c} in \cite{Adamovic}. Below we discuss such a
realization.

Set
\[
 \alpha = Q + \frac{N-1}{N} A_1 + \dots + \frac{1}{N} A_{N-1}, \qquad \beta = Y - \alpha,
\]
elements of the vertex algebra $\calV^0$. Note that we have
$(\alpha, \alpha) = 1$, $(\beta, \beta) = -1$, $(\alpha, \beta) = 0$.
We consider the fermionic vertex operators $\Psi^{\pm}(z) := e^{\pm \intalpha}(z)$
and $e^{\pm \intbeta}(z)$. Note that $H = - \beta$ and $e^{- \XI} = \Psi^{-}_{(-1)} e^{- \intbeta}$ by
definition, and hence we have 
$(\partial - H_{(-1)})^n\, e^{- \XI} = (\partial^n \Psi^{-})_{(-1)} e^{- \beta}$
for $n \ge 0$. Then, we have
\begin{multline*}
 F = - \R{N} \cdots \R{1} e^{- \XI} \\
 = (-1)^{N+1} (\partial - H_{(-1)} + \Mu_{N (-1)}) \cdots (\partial - H_{(-1)} + \Mu_{1 (-1)}) e^{- \XI} \\
 = (-1)^{N+1} \sum_{m=0}^{N} \W_{m (-1)} (\partial - H_{(-1)})^{N-m} e^{- \XI} \\
 = (-1)^{N+1} \sum_{m=0}^{N} \W_{m (-1)} (\partial^{N-m} \Psi^{-})_{(-1)} e^{- \beta}
\end{multline*}
since $H_{(-1)}$ commutes with $\Mu_{i (-1)}$ for $i=1$, $\dots$, $N$.
Therefore, we obtain the following realization of the vertex algebra $\calW^{-N}(\frsl_N, f_{sub})$.

\begin{theorem}[Adamovi\'{c}'s conjecture \cite{Adamovic}]
\label{thm:Adamovic}
 For $N \ge 2$, the subregular $\calW$-algebra $\calW^{-N}(\frsl_N, f_{sub})$ of type $A_{N-1}$ is
 isomorphic to the vertex algebra strongly generated by the following fields:
 \begin{gather*}
  \NO e^{\intbeta}(z) \Psi^{+}(z) \NO,\quad H(z),\quad \W_m(z) \quad (m=2, \dots, N-1), \\
  \text{and } \sum_{m=0}^{N} \NO \W_{m}(z) e^{- \intbeta}(z) \partial^{N-m} \Psi^{-}(z)\NO.
 \end{gather*}
\end{theorem}

\bibliographystyle{halpha} 
\bibliography{refs}

\begin{thebibliography}{KRW03}

\bibitem[ACL18]{ACL2}
Tomoyuki Arakawa, Thomas Creutzig, and Andrew~R. Linshaw.
\newblock W-algebras as coset vertex algebras.
\newblock {\em arXiv preprint}, arXiv:1801.03822:1--38, 2018.

\bibitem[Ada15]{Adamovic}
Dra\v{z}en Adamovi\'{c}.
\newblock Realization of ${W}^{(2)}_n$ algebra and its {W}hittaker modules at
  the critical level.
\newblock {\em private notes}, pages 1--4, 2015.

\bibitem[AM17]{AM}
Tomoyuki Arakawa and Alexander Molev.
\newblock Explicit generators in rectangular affine {$\calW$}-algebras of type
  {$A$}.
\newblock {\em Lett. Math. Phys.}, 107(1):47--59, 2017.

\bibitem[Ara12]{Arakawa}
Tomoyuki Arakawa.
\newblock {$W$}-algebras at the critical level.
\newblock In {\em Algebraic groups and quantum groups}, volume 565 of {\em
  Contemp. Math.}, pages 1--13. Amer. Math. Soc., Providence, RI, 2012.

\bibitem[Ara17]{Arakawa2}
Tomoyuki Arakawa.
\newblock Introduction to {W}-algebras and their representation theory.
\newblock In {\em Perspectives in {L}ie theory}, volume~19 of {\em Springer
  INdAM Ser.}, pages 179--250. Springer, Cham, 2017.

\bibitem[DSK06]{DK}
Alberto De~Sole and Victor~G. Kac.
\newblock Finite vs affine {$W$}-algebras.
\newblock {\em Jpn. J. Math.}, 1(1):137--261, 2006.

\bibitem[FBZ04]{FB}
Edward Frenkel and David Ben-Zvi.
\newblock {\em Vertex algebras and algebraic curves}, volume~88 of {\em
  Mathematical Surveys and Monographs}.
\newblock American Mathematical Society, Providence, RI, second edition, 2004.

\bibitem[FF90]{FF}
Boris Feigin and Edward Frenkel.
\newblock Quantization of the {D}rinfel'd-{S}okolov reduction.
\newblock {\em Phys. Lett. B}, 246(1-2):75--81, 1990.

\bibitem[FL90]{FL}
Vladimir~A. Fateev and Sergei~L. Lukyanov.
\newblock Additional symmetries and exactly solvable models of two-dimensional
  conformal field theory.
\newblock {\em Sov. Sci. Rev. A. Phys.}, 15(1):1--117, 1990.

\bibitem[FMS86]{FMS}
Daniel Friedan, Emil Martinec, and Stephen Shenker.
\newblock Conformal invariance, supersymmetry and string theory.
\newblock {\em Nuclear Phys. B}, 271(1):93--165, 1986.

\bibitem[Fre05]{Frenkel}
Edward Frenkel.
\newblock Wakimoto modules, opers and the center at the critical level.
\newblock {\em Adv. Math.}, 195(2):297--404, 2005.

\bibitem[FS04]{FS}
B.~L. Feigin and A.~M. Semikhatov.
\newblock {$\scW^{(2)}_n$} algebras.
\newblock {\em Nuclear Phys. B}, 698(3):409--449, 2004.

\bibitem[FZ92]{FZ}
Igor~B. Frenkel and Yongchang Zhu.
\newblock Vertex operator algebras associated to representations of affine and
  {V}irasoro algebras.
\newblock {\em Duke Math. J.}, 66(1):123--168, 1992.

\bibitem[Gen17]{Genra1}
Naoki Genra.
\newblock Screening operators for {$\calW$}-algebras.
\newblock {\em Selecta Math. (N.S.)}, 23(3):2157--2202, 2017.

\bibitem[Gen18]{Genra2}
Naoki Genra.
\newblock Screening operators and parabolic inductions for affine w-algebras.
\newblock {\em arXiv preprint}, arXiv:1806.04417:1--49, 2018.

\bibitem[Kac98]{Kac}
Victor Kac.
\newblock {\em Vertex algebras for beginners}, volume~10 of {\em University
  Lecture Series}.
\newblock American Mathematical Society, Providence, RI, second edition, 1998.

\bibitem[KRW03]{KRW}
Victor Kac, Shi-Shyr Roan, and Minoru Wakimoto.
\newblock Quantum reduction for affine superalgebras.
\newblock {\em Comm. Math. Phys.}, 241(2-3):307--342, 2003.

\bibitem[KW04]{KW1}
Victor~G. Kac and Minoru Wakimoto.
\newblock Quantum reduction and representation theory of superconformal
  algebras.
\newblock {\em Adv. Math.}, 185(2):400--458, 2004.

\bibitem[Mol18]{Molev}
Alexander Molev.
\newblock {\em Sugawara operators for classical {L}ie algebras}, volume 229 of
  {\em Mathematical Surveys and Monographs}.
\newblock American Mathematical Society, Providence, RI, 2018.

\bibitem[Pre02]{Premet}
Alexander Premet.
\newblock Special transverse slices and their enveloping algebras.
\newblock {\em Adv. Math.}, 170(1):1--55, 2002.
\newblock With an appendix by Serge Skryabin.

\bibitem[Smi90]{Smith}
S.~P. Smith.
\newblock A class of algebras similar to the enveloping algebra of {${\rm
  sl}(2)$}.
\newblock {\em Trans. Amer. Math. Soc.}, 322(1):285--314, 1990.

\bibitem[Wak86]{Wakimoto}
Minoru Wakimoto.
\newblock Fock representations of the affine {L}ie algebra {$A^{(1)}_1$}.
\newblock {\em Comm. Math. Phys.}, 104(4):605--609, 1986.

\bibitem[Zhu96]{Zhu}
Yongchang Zhu.
\newblock Modular invariance of characters of vertex operator algebras.
\newblock {\em J. Amer. Math. Soc.}, 9(1):237--302, 1996.

\end{thebibliography}

\vspace*{7mm}
\footnotesize{
N.G.: Research Institute for Mathematical Sciences, Kyoto University, Kyoto 606-8502 JAPAN, \\
{\em E-mail address}: \texttt{gnr@kurims.kyoto-u.ac.jp} 

T.K.: Department of Mathematics, Faculty of Pure and Applied Sciences, University of Tsukuba,
Tsukuba, Ibaraki 305-8571, JAPAN. \\
{\em E-mail address}: \texttt{kuwabara@math.tsukuba.ac.jp}
}

\end{document}